\documentclass[graybox]{svmult}

\usepackage{mathptmx}       
\usepackage{helvet}         
\usepackage{courier}        
\usepackage{type1cm}        
\usepackage{makeidx}         
\usepackage{graphicx}        
\usepackage{multicol}        
\usepackage[bottom]{footmisc}
\usepackage{amsfonts}
\usepackage{dsfont}
\usepackage{amsmath}
\usepackage{amssymb}
\usepackage{mathtools}
\usepackage{centernot}
\usepackage{enumerate}
\usepackage{enumitem}
\usepackage[titletoc,title]{appendix}
\usepackage{booktabs,multirow}
\usepackage{natbib}
  \def\citeapos#1{\citeauthor{#1}'s (\citeyear{#1})}
  
\makeindex             

\newcommand{\Var}{\operatorname{\mathsf{Var}}}
\newcommand{\Cov}{\operatorname{\mathsf{Cov}}}
\newcommand{\Corr}{\operatorname{\mathsf{Corr}}}
\newcommand{\MSE}{\operatorname{\mathsf{MSE}}}
\newcommand{\bias}{\operatorname{\mathsf{bias}}}
\newcommand{\CZ}{\mathcal{Z}}
\renewcommand{\E}{\mathds{E}}
\newcommand{\RR}{\mathds{R}}

\newcommand{\os}{{\rm{os}}}
\newcommand{\muos}{\mu_{\os}}
\newcommand{\hmuos}{\h\mu_{\os}}
\newcommand{\bn}{\bm{n}}
\newcommand{\sfc}{\mathsf{C}}
\newcommand{\sfv}{\mathsf{V}}
\newcommand{\ej}{\mathsf{e}_{j}}
\newcommand{\ejl}{\mathsf{e}_{j;l}}
\newcommand{\LL}{\ell}
\newcommand{\hd}{\h{d}}

\newcommand{\h}[1]{\widehat{#1}}

\newcommand{\bm}[1]{\pmb{#1}}
\newenvironment{pr}[1]{\noindent$\textit{#1}$.\ }{\smartqed\bigskip\qed}
\newcommand{\pto}{\xrightarrow{\textrm{~p~}}}

\begin{document}

\title*{The out-of-source error in multi-source cross validation-type procedures}
\author{Georgios Afendras \and Marianthi Markatou}
\authorrunning{Afendras and Markatou} 
\institute{Georgios Afendras \at Department of Biostatistics and Jacobs School of Medicine and Biomedical Sciences, University at Buffalo, Buffalo, NY, USA, \email{gafendra@buffalo.edu}
\and Marianthi Markatou \at Department of Biostatistics and Jacobs School of Medicine and Biomedical Sciences, University at Buffalo, Buffalo, NY, USA, \email{markatou@buffalo.edu}}
%
%
\maketitle


\abstract{
A scientific phenomenon under study may often be manifested by data arising from processes, i.e. sources, that may describe this phenomenon. In this contex of multi-source data, we define the ``out-of-source'' error, that is the error committed when a new observation of unknown source origin is allocated to one of the sources using a rule that is trained on the known labeled data. We present an unbiased estimator of this error, and discuss its variance. We derive natural and easily verifiable assumptions under which the consistency of our estimator is guaranteed for a broad class of loss functions and data distributions. Finally, we evaluate our theoretical results via a simulation study.
}

\section{Introduction}
\label{sec.intr}

In many situations data arise not from a single source but from multiple sources, each of which may have a specific generating process. An example of such a situation is the monitoring and diagnosis of cardiac arrhythmias.

Monitoring devices in cardiac intensive care units use data from electrocardiogram (ECG) channels to automatically diagnose cardiac arrhythmias. However, data from other sources like arterial pressure, ventilation, etc. are often available, and each of these sources has a specific data generating process. Other potential data sources include nuclear medicine tests and echocardiograms. Other examples arise in natural language processing where labeled data for information extraction or parsing are obtained from a limited set of document types.

Cross validation is a fundamental statistical method used extensively in both, statistics and machine learning. A fundamental assumption in using cross validation is that the observations are realizations of exchangeable random variables, that is both, the training and test data come from the same source. However, this is not necessarily the case when data come from multiple sources. \citet{GS2013} state that ``for data of this nature, a common procedure is to arrange the cross validation procedure by source''.

In this setting, we are interested in estimating the generalization error of learning algorithms.

The generalization error is defined as the error an algorithm makes on cases that the algorithm has never seen before, and is important because it relates to the algorithm's prediction capabilities on independent data. The literature includes both, theoretical investigations of risk performance of machine learning algorithms as well as numerical comparisons.

Estimation of the generalization error can be achieved via the use of resampling techniques. The process consists of splitting the available data into a learning or training set and a test set a large number of times and averaging over these repetitions. A very popular resampling technique is cross validation. We are interested in investigating the use of cross validation in the case of multi-source data, where testing occurs on elements that may not have been part of the training set on which the learning algorithm was trained. We do not offer here a detailed overview of cross validation. The interested reader is referred to \citet{Stone1974,Stone1977} for foundational aspects of cross validation, \citet[Ch.s~3,8]{BFOS1984}, \citet{Geisser1975}, and to \citet{AC2010} for a comprehensive survey.

Two very popular cross validation procedures are the $k$-fold cross validation and the random cross validation. In $k$-fold cross validation we split the data randomly in $k$ equal parts. Each of these parts is a test set and the remaining data serves as a training set. In random cross validation case we split randomly the data into two sets of sizes $n_1$ and $n_2$ ($n=n_1+n_2$ is the total sample size) that serve as training/test sets and repeat this process $J$ times.

Carrying out inference on the generalization error requires insight into the variance of the estimator of it. \citet{NB2003} provided estimators for the variance of the random cross validation estimator of the generalization error, while \citet{BG2004} addressed estimators of variance in $k$-fold cross validation. Furthermore, \citet{MTBH2005} proposed moment approximation-based estimators of the same cross validation estimators of the generalization error, and compared these estimators with those provided by \citeauthor{NB2003}.

\citet{AM2016} study the optimality of the data splitting for both, random and $k$-fold cross validation procedures establishing the corresponding optimization rules. Their work offers closed form solutions to corresponding optimization problems for a broad and commonly used class of loss functions.

A recent article by \citet{GS2013} addresses a formulation of the {\it out-of-source} (OOS) error in a multi-source data setting. In their framework there are $k$ sources. If the data size is $n$, the sample size of each source is  $n/k$. The observations of each source are independent and identically distributed (iid) realizations of random variables/vectors from an unknown distribution. The elements that belong to a specific source constitute a test set, while the union of the elements of the remaining sources constitutes the corresponding training set. \citet{GS2013} construct their cross validation-type decision rule using the elements of the aforementioned training set. In this sense, their procedure can be thought of as $k$-fold cross validation with the fundamental difference that the test set data does not necessarily follow the same distribution as the training data.

Recently, multi-source data analysis has received considerable attention in the literature. \citet{B-DBCKPV2010} study the performance of classifiers trained on source data but tested on target data, that is data that do not necessarily follow the same distribution with the source data. Specifically, they study conditions under which a classifier performs well, as well as strategies to combine a small amount of labeled target data at the training step of the classifier to facilitate better performance.

This paper is organized as follows. Section \ref{sec.motiv} presents a motivation of the problems studied here, while Section \ref{sec.F-N} establishes the framework and notation that is used in this article. Section \ref{sec.OOSe} defines the estimator of OOS error and discusses the properties of this estimator. Section \ref{sec.sim} presents simulation results while Section \ref{sec.discussion} offers a discussion. Finally, \ref{append.Isserlis} shows some useful existing results and \ref{append.pr} contains the proofs of the obtained results.

\section{Motivation}
\label{sec.motiv}
\citeauthor{GS2013}'s hypothesis that each source has the same number of elements with each other is too restrictive (and some times not realistic) in practice. Also, the construction of the decision rule based on all of the elements of the training set often leads to various pathologies, as we see in the following example.

\begin{example}
\label{ex.motivation}
Consider a data set of observations that are realizations of independent variables and has size $n=30$, say $\{Z_1,\ldots,Z_{30}\}$. Assume that the data set arises from three sources with $n_1=n_2=n_3=10$ observations each; in general, denote by $n_j$ the sample size associated with the $j$th source. A variable of the first source follows $N(-\mu,1)$, of the second source follows $N(0,1)$ and of the third source follows $N(\mu,1)$. Let the squared error loss be used and suppose that a new variable $Z$ comes from the second source and is independent from the remaining variables in the data set. Additionally, assume that the decision rule is the sample mean. According to \citeapos{GS2013} formulation, $Z$ has an OOS error which is $\E(Z-\overline{Z}_{1,3})^2$, where $\overline{Z}_{1,3}$ is the average of the union of the elements of the first and third sources, that is $\overline{Z}_{1,3}=\frac{1}{20}\sum_{i\in S_1\cup S_3}Z_i$, where $S_j$ denotes the set of indices of the elements of the $j$th source. One can easily see that $Z-\overline{Z}_{1,3}\sim N(0,1.05)$. Therefore, $\E(Z-\overline{Z}_{1,3})^2=\Var(Z-\overline{Z}_{1,3})+\E^2(Z-\overline{Z}_{1,3})=1.05$. Observe that in this case \citeauthor{GS2013}'s formulation has the pathology that the preceding error is independent of the value of $\mu$. We see below, see Example \ref{ex.motivation2}, that the OOS error that is addressed by \citet{GS2013} is significantly different than the actual OOS error.
\end{example}

In view of the above, it is clear that we need to re-formulate the definition of the OOS error in the context of multi-source data to take into account the fact that in practice, not all sources have the same number of observations.

\section{Framework and notation}
\label{sec.F-N}
Assume $k$ sources, where $k$ is a fixed number. Let a data set $\{Z_1,\ldots,Z_n\}$ of size $n$ be observed by the following mechanism: The observations come from the sources that follow a distribution $\bm{p}=(p_1,\ldots,p_k)$. That is, the percentage of observations of the $j$th source is $p_j$, $j=1,\ldots,k$ (the sample size of the $j$th source is $n_j=np_j$). The vector of the numbers of observations of the sources is denoted by $\bm{n}=n\bm{p}=(n_1,\ldots,n_k)$. Each $Z_i$ is labeled by its source and it is independent from the remaining observations whether those come from the same or different sources. The observations of the $j$th source constitute an iid collection of size $n_j$ from an unknown distribution $F_j$.

A new unobserved variable, say $Z$, comes from a source and is independent from $\{Z_1,\ldots,Z_n\}$. The probability of the event ``the variable $Z$ belongs to the $j$th source'' is $p_j=n_j/n$, and follows the distribution $F_j$, $j=1,\ldots,k$. The OOS error is the error that arises between the variable $Z$ and the $k-1$ foreign sources with respect to $Z$, when a loss function $L$ is used for measuring this error.

To formalize the above procedure, first we give a list of definitions and notations. Let $N=\{1,\ldots,n\}$. For each $A\subseteq N$, we denote by $\CZ_A$ the set $\CZ_A\doteq\{Z_i \mid i\in A\}$. The set of the indices of the $j$th source is denoted by $S_j$ and the set of observation of this source is $\CZ_{S_j}$. The loss function $L$ is a measurable nonnegative real function $L(T,\hd)$, where $\hd$ is a decision rule and $T$ is the target variable. The decision rule is constructed using the elements of a set $\CZ_A$, i.e. $\hd=\hd(\CZ_A)$, while the target variable is an element $Z_i\notin\CZ_A$. Hereafter, we write $\hd_{j}\equiv \hd_{j,n}=\hd(\CZ_{S_j})$, $j=1,\ldots,k$, when the decision rule is constructed based on the elements of the $j$th source and the total sample size is $n$.

Since $Z_i$s are independent and the elements of each source are identically distributed, the following are obvious.
\begin{enumerate}[leftmargin=21pt, label=\rm E.\arabic*:,ref=\rm E.\arabic*]
\item \label{eq.exch1}
$\big(L(Z_i,\hd_{l}),L(Z_{i'},\hd_{l})\big)$, $i\in S_j$, $i'\in S_{j'}$ with $j\ne j' \ne l \ne j$, are exchangeable;

\item \label{eq.exch2}
$\big(L(Z_i,\hd_{l}),L(Z_{i},\hd_{l'})\big)$, $i\in S_j$ with $j\ne l \ne l' \ne j$, are exchangeable;

\item \label{eq.exch3}
$\big(L(Z_i,\hd_{l}),L(Z_{i'},\hd_{l'})\big)$, $i\ne i'\in S_j$ with $j\ne l \ne l' \ne j$, are exchangeable;

\item \label{eq.exch4}
$\big(L(Z_i,\hd_{l}),L(Z_{i'},\hd_{j})\big)$, $i\in S_j$, $i'\in S_{j'}$ with $j\ne j' \ne l \ne j$, are exchangeable.

\item \label{eq.exch5}
$L(Z_i,\hd_{l})$ and $L(Z_{i'},\hd_{l'})$, $i\in S_j$, $i'\in S_{j'}$, are independent for all indices $j,j',l,l'$ such that $\{j,l\}\cap\{j',l'\}=\varnothing$.
\end{enumerate}

Now we are in a position to present the algebraic form of the OOS error. Hereafter we assume that the loss function has finite moment of the first order; that is, $\E|L(Z_i,\hd_l)|<\infty$ for all $Z_i\in\CZ_j$ and $l\ne j$. Given that the variable $Z$ comes from the $j$th source and that the decision is constructed based on the elements of the $l$th source, the error committed is
\[
\ejl\doteq\E[L(Z,\hd_{l})] \ \ \textrm{when $Z\sim F_j$ and $Z$, $\CZ_{S_l}$ are independent},
\]
that is, $\ejl$ is the expected value of the loss function when the decision rule is constructed based on the elements of the $l$th source and the target variable belongs to the $j$th source and is independent of the elements of the $l$th source. Taking into account the distribution of the sources and using the conditional total probability theorem, given that the variable $Z$ comes from the $j$th source the error is
\[
\ej\doteq\frac{1}{1-p_j}\sum_{\mathclap{l\ne j}}p_l\ejl;
\]
this is the error that is created from an observation from the $j$th source when compared against observations from the other sources. According to the total probability theorem, the total OOS error is defined by
\begin{equation}
\label{eq.OOSer}
\mu^{(\bn)}_{\os}\doteq\sum_{\mathclap{j=1}}^{k}p_j\ej=\sum_{\mathclap{j=1}}^{k}od_j\sum_{\mathclap{l\ne j}}p_l\ejl,
\end{equation}
where $od_j\doteq p_j/(1-p_j)$ is the odds ratio of the $j$th source. This error can be thought of as a generalization-type error.

\section{The OOS error estimation}
\label{sec.OOSe}

Here, we give an estimator of the OOS error defined in the previews section and investigate the properties of this estimator.

\subsection{Estimating the OOS error}
\label{ssec.OOSest}
We are interested in estimating the OOS error. To simplify notation we use $\LL_{i,j}$ to denote $L(Z_i,\hd_{j})$. By definition, a natural estimator of $\ejl$ is $\frac{1}{n_j}\sum_{i\in S_j}\LL_{i,l}$, and thus, a natural estimator of $\mu^{(\bn)}_{\os}$ is
\begin{equation}
\label{eq.OOSest}
\h\mu^{(\bn)}_{\os}\doteq\frac{1}{n}\sum_{\mathclap{j=1}}^{k}\frac{1}{n-n_j}\sum_{\mathclap{l\ne j}}n_l\sum_{\mathclap{i\in S_j}}\LL_{i,l}.
\end{equation}
This estimator is a cross validation-type estimator of the OOS error. When $k$-fold cross validation is used to estimate generalization error the data are split in $k$ equal parts. Each of these parts is a test set and its complement set is the corresponding training set. The target variable is a variable of the test set and the decision rule is constructed based on the elements of the corresponding training set \citep[for more details see, for example,][]{AM2016}. Here, we have $k$ test sets ($\CZ_{S_1},\ldots,\CZ_{S_k}$), which are defined by the labeling of the data and are a partitioning of the data set. For each test set $\CZ_{S_j}$ the corresponding training set $\CZ_{S_{j}}^c=\bigcup_{l\ne j}\CZ_{S_{l}}$ is partitioned into $k-1$ training sub-sets $\CZ_{S_1},\ldots,\CZ_{S_{j-1}},\CZ_{S_{j+1}},\ldots,\CZ_{S_k}$. The target variable is a variable of the test set and for each $l\ne j$ the decision rule is constructed based on the elements of the training sub-set $\CZ_{S_{l}}$.

Hereafter we write $\muos$ and $\hmuos$ instead $\mu^{(\bn)}_{\os}$ and $\h\mu^{(\bn)}_{\os}$ respectively. The following example illustrates the difference between the OOS error that introduced by \citeauthor{GS2013} and that we have defined in relationship \eqref{eq.OOSer}.

\begin{example}[Example \ref{ex.motivation} continued]
\label{ex.motivation2}
Let the data be as in Example \ref{ex.motivation} and the squared error loss is used. Let us consider $Z^{(1)}\sim N(-\mu,1)$, $Z^{(2)}\sim N(0,1)$, $Z^{(3)}\sim N(\mu,1)$ and $Z^{(1)}$, $Z^{(2)}$, $Z^{(3)}$, $\{Z_1,\ldots,Z_n\}$ are independent. Then, $Z^{(1)}-\overline{Z}_{2,3}\sim N(-3\mu/2,1.05)$, $Z^{(2)}-\overline{Z}_{1,3}\sim N(0,1.05)$ and $Z^{(3)}-\overline{Z}_{1,2}\sim N(3\mu/2,1.05)$. Hence, the OOS error given by \citet{GS2013} is
\[
\mu_{\rm CVS}=\frac{1}{3}\left\{\E\!\big(Z^{(1)}-\overline{Z}_{2,3}\big)^2+\E\!\big(Z^{(2)}-\overline{Z}_{1,3}\big)^2+\E\!\big(Z^{(3)}-\overline{Z}_{1,2}\big)^2\right\}=1.05+\frac{3}{2}\mu^2.
\]
Using the more general Example \ref{exm.normal,sq-abs}\eqref{exm.normal,sq-abs(a)} below, relation \eqref{eq.normal,sq-abs1} gives that the OOS error given by \eqref{eq.OOSer} is $\muos=1.1+2\mu^2$.
\end{example}

\subsection{Bias and variance  of $\bm{\hmuos}$}
\label{ssec.bias,var}

In this section we investigate the bias and variance, and so the mean square error, of the OOS error estimator $\hmuos$. In view of \ref{eq.exch1}--\ref{eq.exch5}, we state and prove the following theorem.

\begin{theorem}
\label{theo.bias,var}
Assume that $\E[L(Z,\hd_l)]^2<\infty$ when $Z\sim F_j$ and $Z$, $\CZ_{S_l}$ are independent for all $j\ne l$. Then,
\begin{enumerate}[label=\rm(\alph*), ,ref=\alph*]
\item
\label{theo.bias,var(a)}
the estimator $\hmuos$ given by \eqref{eq.OOSest} is an unbiased estimator of the OOS error;
\item
\label{theo.bias,var(b)}
the variance of $\hmuos$ is
\[
\begin{split}
\Var(\hmuos)=
&\sum_{\mathclap{j=1}}^k \frac{od_j^2}{np_j}\left(\sum\limits_{\mathclap{l\ne j}} p_l^2\big(\sfv_{j;l}+n(p_j-1/n)\sfc_{j;l}\big)+\mathop{\sum\sum}\limits_{\mathclap{l\ne l' \colon l,l'\ne j}}p_lp_{l'}\sfc_{j;l,l'}\right)\\
&+\mathop{\sum\sum}\limits_{\mathclap{j\ne j'}}od_jod_{j'}\left(\ \ \sum\limits_{\mathclap{l\ne j,j'}}p_l(p_l\sfc_{j,j';l}+2p_j\sfc_{j,j';l,j})\right),
\end{split}
\]
where $\sfv_{j;l}=\Var L(Z,\hd_{l})$ when $Z\sim F_j$ and is independent of $\CZ_{S_l}$; $\sfc_{j;l}=\Cov(L(Z,\hd_{l}),L(Z',\hd_{l}))$ when $Z$ and $Z'$ are iid from $F_j$ and are independent of $\CZ_{S_l}$; $\sfc_{j;l,l'}=\Cov(L(Z,\hd_{l}),L(Z,\hd_{l'}))$ when $Z\sim F_j$ and is independent of $\CZ_{S_l}\cup\CZ_{S_{l'}}$; $\sfc_{j,j';l}=\Cov(L(Z,\hd_{l}),L(Z',\hd_{l}))$ when $Z\sim F_j$, $Z'\sim F_{j'}$ and $Z$, $Z'$, $\CZ_{S_l}$ are independent; and $\sfc_{j,j';l,j}=\Cov(L(Z,\hd_{l}),L(Z',\hd_{j}))$ when $Z\in \CZ_{S_j}$, $Z'\sim F_{j'}$ and $Z'$, $\CZ_{S_j}$, $\CZ_{S_l}$ are independent.
\end{enumerate}
\end{theorem}

Now we investigate the consistency of the estimator $\hmuos$. First, we are interested in finding simple and natural conditions that imply the desired result. Very often the sequence, with respect to $n$, of the decision rules $\hd_{j;n}$ converges in probability to a constant for each $j$. For example, if the decision rule is the sample mean of the elements of the $j$th source, say $\overline{Z}_j$, and $F_j$ has mean $\mu_j$, then $\overline{Z}_j\pto\mu_j$. Also, the finiteness of the variance of the OOS error estimator requires that $\E[L(Z_i,\hd_{l})]^2<\infty$ for all $Z_i\in\CZ_{S_j}$ and $j\ne l$. In view the above observations, we state the following conditions/assumptions:
\begin{enumerate}[leftmargin=21pt, label=\rm C.\arabic*:,ref=\rm C.\arabic*]
\item \label{cond1}
$\hd_{j;n}\pto d_j$, as $n\to\infty$, for all $j=1,\ldots,k$, where $d_j$s are constants.

\item \label{cond2}
There exist $\theta,M>0$ such that $\E[L(Z,\hd_{l,n})]^{2+\theta}\le M$ when $Z\sim F_j$ and $Z$, $\CZ_{S_l}$ are independent, for all $j\ne l$ and $n$.
\end{enumerate}

\begin{theorem}
\label{theo.consistency}
Let $L$ be a continuous loss function and suppose that \ref{cond1}, \ref{cond2} hold. Then, $\Var(\hmuos)\to0$ as $n\to\infty$ and, thus, $\hmuos$ is a consistent estimator of $\muos$.
\end{theorem}

The following Examples \ref{exm.abs,sq} and \ref{exm.normal,sq-abs} show the usefulness of Theorem \ref{theo.consistency}.

\begin{example}
\label{exm.abs,sq}
Let $\CZ_{S_j}$ be an iid collection of random variables (rv's) from $F_j$, $j=1,\ldots,k$, the decision rules are the usual averages of the elements of the sources, $F_j$ does not depend on $n$ and has mean $\mu_j$ and variance $\sigma_j^2$.

\begin{enumerate}[leftmargin=14pt, topsep=0pt, label=\rm(\alph*),ref=\rm\alph*]
\item \label{exm.abs,sq(a)}
Let the absolute error loss be used. Suppose that $F_j$ has finite moments of order $2+\theta_j$ for some $\theta_j>0$, $j=1,\ldots,k$. Then, $\hd_j=\overline{Z}_j=\frac{1}{n_j}\sum_{i\in S_j}Z_i\pto \mu_j$, that is, \ref{cond1} is satisfied. Set $\theta=\min_{j=1,\ldots,k}{\{\theta_j\}}>0$, $\beta_{2+\theta}=\max_{j=1,\ldots,k}\{\E|Z|^{2+\theta}$  when $Z\sim F_j\}<\infty$ and $M=2^{3+\theta}\beta_{2+\theta}<\infty$. For each $j\ne l$ and $Z\sim F_j$ such that $Z$, $\CZ_{S_l}$ are independent we have that $\E|L(Z,\hd_{l,n})|^{2+\theta}=\E|Z-\overline{Z}_l|^{2+\theta}\le 2^{2+\theta}\left(\E|Z|^{2+\theta}+\E|\overline{Z}_l|^{2+\theta}\right)\le 2^{2+\theta}\left(\E|Z|^{2+\theta}+\frac{1}{n_l}\sum_{i\in S_l}\E|Z_i|^{2+\theta}\right)\le2^{3+\theta}\beta_{2+\theta}=M$; and so \ref{cond2} is satisfied. Therefore, Theorems \ref{theo.bias,var} and \ref{theo.consistency} show that $\MSE(\hmuos)=\bias^2(\hmuos)+\Var(\hmuos)=\Var(\hmuos)\to0$ as $n\to\infty$.

\item \label{exm.abs,sq(b)}
Let the squared error loss be used, that is $L(Z,\hd_{l,n})=(Z-\overline{Z}_l)^2$, and suppose that $F_j$ has finite moments of order $4+\theta_j$ for some $\theta_j>0$, $j=1,\ldots,k$. Using the same arguments as in \eqref{exm.abs,sq(a)}, we obtain that $\MSE(\hmuos)\to0$ as $n\to\infty$.
\end{enumerate}
\end{example}

\begin{example}
\label{exm.normal,sq-abs}
Suppose $F_j\sim N(\mu_j,\sigma_j^2)$ and $\hd_j=\overline{Z}_j=\frac{1}{n_j}\sum_{i\in S_j}Z_i$.
\begin{enumerate}[leftmargin=14pt, topsep=0pt, label=\rm(\alph*),ref=\rm\alph*]
\item \label{exm.normal,sq-abs(a)}
Let the squared error loss be used. Then, we calculate (see in \ref{append.pr})
\begin{equation}
\label{eq.normal,sq-abs1}
\muos=\sum_{\mathclap{j=1}}^k p_j\sigma_j^2 + \sum_{\mathclap{j=1}}^k od_j\sum_{\mathclap{l\ne j}}p_l(\mu_j-\mu_l)^2 + \frac{1}{n}\sum_{\mathclap{j=1}}^k od_j\sum_{\mathclap{l\ne j}}\sigma_j^2;
\end{equation}
and the quantities $\sfv_{j;l}$, $\sfc_{j;l}$, $\sfc_{j;l,l'}$ $\sfc_{j,j';l}$ and $\sfc_{j,j';l,j}$ that appear in the variance of $\hmuos$ in Theorem \ref{theo.bias,var}\eqref{theo.bias,var(b)} are
\begin{subequations}
\label{eq.normal,sq-abs2}
\begin{equation}
\label{eq.normal,sq-abs2(a)}
\sfv_{j;l}=2\left(\sigma_j+\frac{\sigma_l^2}{np_l}\right)\left[\sigma_j+\frac{\sigma_l^2}{np_l}+2(\mu_j-\mu_l)^2\right],
\end{equation}
\begin{equation}
\label{eq.normal,sq-abs2(b)}
\left.
\begin{split}
&\sfc_{j;l}=2\frac{\sigma_l^2}{np_l}\left(\frac{\sigma_l^2}{np_l}+2(\mu_j-\mu_l)^2\right),\\
&\sfc_{j;l,l'}=2\sigma_j^2\left(\sigma_j^2+2(\mu_j-\mu_l)(\mu_j-\mu_{l'})\right),\\
&\sfc_{j,j';l}=2\frac{\sigma_l^2}{np_l}\left(\frac{\sigma_l^2}{np_l}+2(\mu_j-\mu_l)(\mu_{j'}-\mu_l)\right),\\
&\sfc_{j;l,l'}=2\frac{\sigma_j^2}{np_j}\left(\frac{\sigma_j^2}{np_j}-2(\mu_j-\mu_l)(\mu_{j'}-\mu_j)\right).\\
\end{split}
\right\}
\end{equation}
\end{subequations}
Observe that $\sfc_{j;l},\sfc_{j,j';l},\sfc_{j,j';l,j}\to 0$ as $n\to\infty$; specifically, these covariances are $O(1/n)$ functions as $n\to\infty$. It is obvious that $\Var(\hmuos)=O(1/n)$ as $n\to\infty$. This example is a confirmation of Theorem \ref{theo.consistency} for this case.
\item
\label{exm.normal,sq-abs(b)}
Let the absolute error loss be used. Then, see in \ref{append.pr},
\begin{equation}
\label{eq.normal,sq-abs3}
\muos=\sum_{\mathclap{j=1}}^{k}od_j\sum_{\mathclap{l\ne j}}p_l \left\{\mu_{j;l}\left[1-2\Phi\left(-\frac{\mu_{j;l}}{\sigma_{j;l}}\right)\right]+\sigma_{j;l}\sqrt{\frac2\pi}\exp\left(-\frac{\mu_{j;l}^2}{2\sigma_{j;l}^2}\right)\right\},
\end{equation}
where $\mu_{j;l}=\mu_j-\mu_l$, $\sigma_{j;l}^2=\sigma_j^2+\sigma_l^2/(np_l)$ and $\Phi$ denotes the cumulative distribution function of the standard normal distribution. The calculations of the covariances $\sfc_{j;l}$, $\sfc_{j;l,l'}$, $\sfc_{j,j';l}$ and $\sfc_{j;l,l'}$ in Theorem \ref{theo.bias,var}(b) are rather difficult.
\end{enumerate}
\end{example}

In practice, the data distributions of the sources are unknown and, thus, the OOS error must be estimated. They are loss function for which the calculation of the variance of $\hmuos$ is difficult, or impossible, even if the distribution of the data sources is known, cf. Example \ref{eq.normal,sq-abs1}\eqref{exm.normal,sq-abs(b)}. Furthermore, in the formulation of Example \ref{eq.normal,sq-abs1} where the absolute error loss is used, if we consider $F_j\sim U(a_j,b_j)$, $j=1,2,3$, the calculation of OOS error in closed form is impossible.

\subsection{On variance estimation}
\label{ssec.VarEst}

Here we investigate the possibility of the variance estimation of the OOS error.

First, we present a general result and some useful observations that arise from it. \citet{NB2003} study the variance estimation of the random cross validation estimator of the generalization error of a computer algorithm when $L(Z_i,\hd_j)$ for all realizations $\CZ_{S_j}$ and $Z_i\in\CZ_{S_j}^c$ are exchangeable. They prove that ``There is no general unbiased estimator of the variance of the random cross validation estimator that involves the $L(Z_i,\hd_j)$s in a quadratic and/or linear way.'' \citep[see][Proposition~3, page~246]{NB2003}. This result holds in a more general form.

\begin{lemma}
\label{lem.moments}
Let $X_1,X_2,\ldots,X_n$ be a collection of random variables. If $\E(X_j)=\mu$, $\Var(X_j)=\sigma^2$ and $\Cov(X_j,X_{j'})=\sfc$, $j\ne j'$, are unknown parameters, then we can find unbiased estimators of the second moments of $X_j$s only for the cases of linear combinations of $\sigma^2+\mu^2$, $\sfc+\mu^2$.
\end{lemma}

The following corollary follows immediately from Lemma \ref{lem.moments}.

\begin{corollary}
\label{cor.moments}
Let $X_1,X_2,\ldots,X_n$ be as in Lemma \ref{lem.moments}. Then, {\rm(a)} there does not exist an unbiased estimator of $\Var(\overline{X})$, where $\overline{X}$ is the usual average of $X_j$s, and {\rm(b)} there do not exist unbiased estimators of $\mu^2$, $\sigma^2$, $\sfc$.
\end{corollary}

\begin{remark}
\label{rem.moments}
(a) If one of the parameters $\mu^2$, $\sigma^2$ and $\sfc$ is known, then we can provide unbiased estimators for each linear combination of the other two parameters.
\medskip

\noindent
(b) The statistic $s^2=\frac{1}{n-1}\sum_{j=1}^{n}(X_j-\overline{X})^2$ is an unbiased estimator of $\sigma^2-\sfc$. The variance of this estimator is
\[
\Var(s^2)=\frac{1}{(n-1)^2}\left\{\sum_{\mathclap{j=1}}^{n}\Var(X_j-\overline{X})^2+\mathop{\sum\sum}\limits_{\mathclap{1\le j< j'\le n}}\Cov\big((X_j-\overline{X})^2,(X_{j'}-\overline{X})^2\big)\right\}.
\]
It is possible $\Var(s^2) \centernot\to0$, as $n\to\infty$, see Example \ref{ex.Y,e} below; and thus, $s^2$ is not consistent estimator of $\sigma^2-\sfc$.
\medskip

\noindent
(c) In both random and $k$-fold cross validation estimators of the generalization error of a computer algorithm, the sequence of the test set errors are as in Lemma \ref{lem.moments} \citep[see][Proposition 1]{AM2016}. For both of these cases the cross validation estimator is the usual average of the test set errors. Thus, the unbiased estimation of the variance of the cross validation estimator is impossible. Let $\h{\mu}_j$, $j=1,\ldots,J$, denote the test set errors and $\h{\mu}_{{\rm{CV},J}}=\frac{1}{J}\sum_{j=1}^{J}\h{\mu}_j$ be the cross validation estimator in the random cross validation procedure. If $s^2_{\h{\mu}_j}=\frac{1}{J-1}\sum_{j=1}^{J}\left(\h{\mu}_j-\h{\mu}_{{\rm{CV},J}}\right)$ and $\rho=\Corr(\h{\mu}_j,\h{\mu}_{j'})$, \citet[p.~248]{NB2003} state that ``\ldots $\left(1+\frac{\rho}{1-\rho}\right)s^2_{\h{\mu}_j}$ is an unbiased estimator of $\Var(\h{\mu}_{{\rm{CV},J}})$''; this sentence is incorrect because the parameter $\rho$ is unknown and, thus, $\left(1+\frac{\rho}{1-\rho}\right)s^2_{\h{\mu}_j}$ is not an estimator (statistic). Of course, it is a random variable with $\E\left[\left(1+\frac{\rho}{1-\rho}\right)s^2_{\h{\mu}_j}\right]=\Var(\h{\mu}_{{\rm{CV},J}})$. In general, the estimation of the correlation $\rho$ is difficult. Nevertheless, even in the case in which we find an unbiased estimator of $\rho$, say $\h{\rho}$, then $\left(1+\frac{\h\rho}{1-\h\rho}\right)s^2_{\h{\mu}_j}$ is not an unbiased estimator of $\Var(\h{\mu}_{{\rm{CV},J}})$, except if $\left(1+\frac{\h\rho}{1-\h\rho}\right)$ and $s^2_{\h{\mu}_j}$ are uncorrelated. Moreover, if $\h\rho$ is consistent estimator of $\rho$, then $\left(1+\frac{\h\rho}{1-\h\rho}\right)s^2_{\h{\mu}_j}$ might is not a consistent estimator of $\Var(\h{\mu}_{{\rm{CV},J}})$, cf.\ (b).
\medskip

\noindent
(d) \citet{MTBH2005} provide moment approximation estimators for the variance of the test set errors and for their covariance in a broad and often used class of cross validation procedures, in both random and $k$-fold cross validation cases. In view of (c), it is clear that their results are very important in practice.
\end{remark}

\begin{example}
\label{ex.Y,e}
Let $0<\sfc<\sigma^2<\infty$ and $\mu\in\RR$. Assume that $Y_1,\ldots,Y_n$ is an iid collection from the distribution with probability mass function $p_Y\left(-(\sigma^2-\sfc)^{1/2}/2\right)=p_Y\left((\sigma^2-\sfc)^{1/2}/2\right)=(n^2-1)/(2n^2-1/2)$, $p_Y\left(-n(\sigma^2-\sfc)\right)=p_Y\left(n(\sigma^2-\sfc)\right)=3/(8n^2-2)$; and $\epsilon\sim N(\mu,\sfc)$ which is independent to $Y_j$s. By straightforward calculations, $\E(Y)=0$, $\Var(Y)=\sigma^2-\sfc$. Consider the rv's $X_j=Y_j+\epsilon$, $j=1,\ldots,n$. Then, one can easily see that the $X_j$s are exchangeable with $\E(X_j)=\mu$, $\Var(X_j)=\sigma^2$ and $\Cov(X_j,X_{j'})=\sfc$ for all $j\ne j'$. By definition of the $X_j$s, $s^2_X=\frac{1}{n-1}\sum_{j=1}^n(X_j-\overline{X})=\frac{1}{n-1}\sum_{j=1}^n(Y_j-\overline{Y})=s^2_Y$. So, $\Var(s^2_X)=\Var(s^2_Y)=\frac{\mu_4^{(Y)}}{n}+\frac{(n-3)\Var(Y)}{n(n-1)}=\frac{12n^4+n^2-1}{n(16n^2-4)}+\frac{(n-3)(\sigma^2-\sfc)}{n(n-1)}\to\infty$ as $n\to\infty$.
\end{example}

In view of Theorem \ref{theo.bias,var}\eqref{theo.bias,var(a)}, Lemma \ref{lem.moments} and Corollary \ref{cor.moments}, the unbiased estimation of the variance of $\hmuos$ is impossible because the quantities $\sfv_{j;l}$, $\sfc_{j;l}$, $\sfc_{j;l,l'}$, $\sfc_{j,j';l}$ and $\sfc_{j,j';l,j}$ are as in Lemma \ref{lem.moments}. For example, let the $j$th source and the $l$th source be two deferent sources. Then, $\{\LL_{i,l}, \ i\in S_j\}$ is a set of exchangeable rv's of size $n_j$ with unknown mean, say $\mu_{j,l}$, variance $\sfv_{j,l}$ and covariance between two elements $\sfc_{j,l}$.

It is a fact that there does not exist a general unbiased estimator of the variance of the OOS error. If someone needs an estimator of the variance of the OOS error for some reason (for example, for statistical inference on the OOS error), one may resort to the bootstrap resampling technique or can follow the moment approximation method of \citet{MTBH2005} when it is possible. Notice that the bootstrap resampling technique in this formulation has a very large computational cost.

\section{Simulation study}
\label{sec.sim}

 Assume we have $k=3$ sources with probability vector $\bm{p}=(0.2,0.3,0.5)$, and thus odds vector $\bm{od}=(1/4,3/7,1)$. Suppose that the elements of each source are iid rv's from a distribution and the squared or absolute error loss is used.

Table \ref{tab.norm} presents the true value of the OOS error, the empirical mean and the empirical mean squared error of the OOS error estimator $\hmuos$, when $N=10^4$ Monte Carlo (M-C) repetitions are used, for various values of the sample size $n$. The elements of each source are normally distributed and the squared and absolute error loss are used. In this case we have the explicit expressions of $\muos$ given by the relations \eqref{eq.normal,sq-abs1} and \eqref{eq.normal,sq-abs3} for both cases of squared and absolute error loss respectively. We observe that for both cases, squared and absolute error loss, the empirical mean square error of $\hmuos$ tends to zero as $n$ tends to infinity, confirming the statements of Theorems \ref{theo.bias,var}(a) and \ref{theo.consistency}.
\begin{table}
\caption{The OOS error, $\muos$, the average of $\hmuos$ and its empirical mean square error $\widehat{\MSE}(\hmuos)$, for $N=10^4$ M-C repetitions, when $k=3$, $\bm{p}=(0.2,0.3,0.5)$, $F_1\sim N(0,9)$, $F_2\sim N(2,1)$, $F_3\sim N(5,5)$ and the squared/absolute error loss is used, for various values of $n$.}
\label{tab.norm}
\begin{tabular*}{\linewidth}
 {@{\hspace{0ex}}l@{\hspace{2ex}}l@{\extracolsep{\fill}}r@{\hspace{0ex}}r@{\hspace{0ex}}r@{\hspace{0ex}}r@{\hspace{0ex}}r@{\hspace{0ex}}r@{\hspace{0ex}}r@{\hspace{0ex}}}
 \addlinespace
 \toprule
 & $n$ & $100$~~ & $200$~~ & $300$~~ & $500$~~ & $700$~~ & $10^3$~~ & $10^4$~~\\
 \hline
 \multirow{5}{*}{\rotatebox{90}{squared}}
 & $\muos$                         &     18.171 &     18.084 &     18.055 &     18.031 &     18.021 &     18.014 &     17.998 \\
 & $\hmuos   $                     &     18.164 &     18.084 &     18.055 &     18.038 &     18.038 &     18.017 &     17.999 \\
 & $\widehat{\bias}\,\!^2(\hmuos)$ & $<10^{-3}$ & $<10^{-3}$ & $<10^{-4}$ & $<10^{-4}$ & $<10^{-5}$ & $<10^{-5}$ & $<10^{-6}$ \\
 & $\widehat{\Var}(\hmuos)$        &     11.524 &      5.897 &      3.873 &      2.314 &      1.664 &      1.177 &      0.117 \\
 & $\widehat{\MSE}(\hmuos)$        &     11.524 &      5.897 &      3.873 &      2.314 &      1.664 &      1.177 &      0.117 \\
 \cline{2-9}
 \multirow{5}{*}{\rotatebox{90}{absolute}}
 & $\muos$                         &      3.639 &      3.636 &      3.635 &      3.635 &      3.634 &      3.634 &      3.633 \\
 & $\hmuos$                        &      3.635 &      3.634 &      3.638 &      3.635 &      3.633 &      3.634 &      3.633 \\
 & $\widehat{\bias}\,\!^2(\hmuos)$ & $<10^{-5}$ & $<10^{-5}$ & $<10^{-5}$ & $<10^{-5}$ & $<10^{-6}$ & $<10^{-7}$ & $<10^{-8}$ \\
 & $\widehat{\Var}(\hmuos)$        &      0.148 &      0.074 &      0.050 &      0.030 &      0.021 &      0.014 &      0.001 \\
 & $\widehat{\MSE}(\hmuos)$        &      0.148 &      0.074 &      0.050 &      0.030 &      0.021 &      0.014 &      0.001 \\
 \bottomrule
 \end{tabular*}
\end{table}

Tables \ref{tab.U}--\ref{tab.G} present the empirical mean and the empirical variance of the OOS error estimator $\hmuos$, for $N=10^4$ M-C repetitions for various values of the sample size $n$, when the elements of each source are uniformly distributed (Table \ref{tab.U}), Student distributed (Table \ref{tab.t}) and gamma distributed (Table \ref{tab.G}), and the squared and absolute error loss are used. For both cases of loss function, squared and absolute error loss, and for all cases of the sources' distributions the empirical empirical variance of $\hmuos$ tends to zero as $n$ tends to infinity, confirming empirically the statement of Theorem \ref{theo.consistency}. Note that for these cases of the sources' distributions we do not have explicit forms of $\muos$ and thus, we cannot present the values $\widehat{\bias}\,\!^2(\hmuos)$ and $\widehat{\MSE}(\hmuos)$. On  the other hand, since $\hmuos$ is an unbiased estimator of $\muos$, for large values of $n$ we have that $\widehat{\MSE}(\hmuos)\approx \widehat{\Var}(\hmuos)$.

\begin{table}
\caption{The average of $\hmuos$ and its empirical variance $\widehat{\Var}(\hmuos)$, for $N=10^4$ M-C repetitions, when $k=3$, $\bm{p}=(0.2,0.3,0.5)$, $F_1\sim U(-1,1)$, $F_2\sim U(1/2,3/2)$, $F_3\sim U(3,7)$ and the squared/absolute error loss is used, for various values of $n$.}
\label{tab.U}
\begin{tabular*}{\linewidth}
 {@{\hspace{0ex}}l@{\hspace{2ex}}l@{\extracolsep{\fill}}r@{\hspace{0ex}}r@{\hspace{0ex}}r@{\hspace{0ex}}r@{\hspace{0ex}}r@{\hspace{0ex}}r@{\hspace{0ex}}r@{\hspace{0ex}}}
 \addlinespace
 \toprule
 & $n$ & $100$~~ & $200$~~ & $300$~~ & $500$~~ & $700$~~ & $10^3$~~ & $10^4$~~\\
 \hline
 \multirow{2}{*}{\rotatebox{30}{squared}}
 & $\hmuos$                 & 17.294 & 17.266 & 17.277 & 17.271 & 17.273 & 17.278 & 17.275 \\
 & $\widehat{\Var}(\hmuos)$ &  1.699 &  0.871 &  0.568 &  0.332 &  0.243 &  0.173 &  0.017 \\
 \cline{2-9}
 \multirow{2}{*}{\rotatebox{30}{absolute}}
 & $\hmuos$                 & 3.8435 & 3.8451 & 3.8425 & 3.8433 & 3.8433 & 3.8430 & 3.8428 \\
 & $\widehat{\Var}(\hmuos)$ & 0.0229 & 0.0113 & 0.0077 & 0.0047 & 0.0033 & 0.0022 & 0.0002 \\
 \bottomrule
 \end{tabular*}
\end{table}

\begin{table}
\caption{The average of $\hmuos$ and its empirical variance $\widehat{\Var}(\hmuos)$, for $N=10^4$ M-C repetitions, when $k=3$, $\bm{p}=(0.2,0.3,0.5)$, $F_1\sim t_7$, $F_2\sim t_5(2)$, $F_3\sim t_6(5)$ (where $t_\nu(\mu)\stackrel{\rm d}{=} t_\nu+\mu$) and the squared/absolute error loss is used, for various values of $n$.}
\label{tab.t}
\begin{tabular*}{\linewidth}
 {@{\hspace{0ex}}l@{\hspace{2ex}}l@{\extracolsep{\fill}}r@{\hspace{0ex}}r@{\hspace{0ex}}r@{\hspace{0ex}}r@{\hspace{0ex}}r@{\hspace{0ex}}r@{\hspace{0ex}}r@{\hspace{0ex}}}
 \addlinespace
 \toprule
 & $n$ & $100$~~ & $200$~~ & $300$~~ & $500$~~ & $700$~~ & $10^3$~~ & $10^4$~~\\
 \hline
 \multirow{2}{*}{\rotatebox{30}{squared}}
 & $\hmuos$                 & 14.976 & 14.927 & 14.954 & 14.943 & 14.936 & 14.936 & 14.925 \\
 & $\widehat{\Var}(\hmuos)$ &  2.677 &  1.376 &  0.941 &  0.548 &  0.388 &  0.274 &  0.027 \\
 \cline{2-9}
 \multirow{2}{*}{\rotatebox{30}{absolute}}
 & $\hmuos$                 & 3.5182 & 3.5142 & 3.5153 & 3.5152 & 3.5150 & 3.5153 & 3.5147 \\
 & $\widehat{\Var}(\hmuos)$ & 0.0417 & 0.0209 & 0.0141 & 0.0084 & 0.0061 & 0.0042 & 0.0004 \\
 \bottomrule
 \end{tabular*}
\end{table}

\begin{table}
\caption{The average of $\hmuos$ and its empirical variance $\widehat{\Var}(\hmuos)$, for $N=10^4$ M-C repetitions, when $k=3$, $\bm{p}=(0.2,0.3,0.5)$, $F_1\sim \exp(1)$, $F_2\sim \Gamma(2,1)$, $F_3\sim \Gamma(10,2)$ and the squared/absolute error loss is used, for various values of $n$.}
\label{tab.G}
\begin{tabular*}{\linewidth}
 {@{\hspace{0ex}}l@{\hspace{2ex}}l@{\extracolsep{\fill}}r@{\hspace{0ex}}r@{\hspace{0ex}}r@{\hspace{0ex}}r@{\hspace{0ex}}r@{\hspace{0ex}}r@{\hspace{0ex}}r@{\hspace{0ex}}}
 \addlinespace
 \toprule
 & $n$ & $100$~~ & $200$~~ & $300$~~ & $500$~~ & $700$~~ & $10^3$~~ & $10^4$~~\\
 \hline
 \multirow{2}{*}{\rotatebox{30}{squared}}
 & $\hmuos$                 & 12.095 & 12.068 & 12.067 & 12.060 & 12.038 & 12.048 & 12.043 \\
 & $\widehat{\Var}(\hmuos)$ &  2.763 &  1.326 &  0.863 &  0.537 &  0.385 &  0.265 &  0.027 \\
 \cline{2-9}
 \multirow{2}{*}{\rotatebox{30}{absolute}}
 & $\hmuos$                 & 3.0669 & 3.0636 & 3.0655 & 3.0662 & 3.0657 & 3.0647 & 3.0652 \\
 & $\widehat{\Var}(\hmuos)$ & 0.0469 & 0.0240 & 0.0163 & 0.0095 & 0.0070 & 0.0048 & 0.0005 \\
 \bottomrule
 \end{tabular*}
\end{table}

\section{Discussion}
\label{sec.discussion}

In this paper we discuss the definition, estimation and properties of the proposed estimator of the out-of-source error in the context of multi-source data, when it is not assumed that all sources have exactly the same number of observations and do not necessarily follow the same distribution. We show that our proposed estimator is unbiased, and we offer natural and easy to verify in practice conditions under which the estimator we propose is consistent.

Most research, both theoretical and empirical, assumes that a learning algorithm is trained and tested using data that follow the same distribution. This setting has been extensively studied in the literature, and uniform convergence theory guarantees that a learning algorithm's empirical error is close to its true error under appropriate assumptions. However, in many practical situations we wish to train a learning algorithm under one or more source domains and then test it on a domain that is potentially different from the source domains. Our work, presented here, studies the out-of-source error in this setting. We further supplement the theoretical results we present here with a simulation that essentially verifies these results.

\begin{acknowledgement}
Dr. Markatou would like to thank the Jacobs School of Medicine and Biomedical Science for facilitating this work through institutional financial resources (to M.\ Markatou) that supported the work of the first author of this paper.
\end{acknowledgement}

\appendix
\renewcommand\thesection{Apendix \Alph{section}}
\renewcommand\theequation{\Alph{section}.\arabic{equation}}
\renewcommand{\thetheorem}{\Alph{section}.\arabic{theorem}}
\renewcommand{\thecorollary}{\Alph{section}.\arabic{corollary}}
\setcounter{equation}{0}
\setcounter{theorem}{0}
\setcounter{corollary}{0}
\begin{appendices}

\section{\ \ On moments of bivariate normal distribution}
\label{append.Isserlis}

\begin{theorem}[\citeauthor{Isserlis1918}, \citeyear{Isserlis1918}]
\label{theo.Isserlis}
Let $(Y_1,Y_2)\sim N_2\left(\bm{0},\left({\sigma_1^2\atop\sigma_{12}}~{\sigma_{12}\atop\sigma_2^2}\right)\right)$. Then, $\E(Y_i^4)=3\sigma_{i}^4$, $\E(Y_1^2Y_2^2)=\sigma_1^2\sigma_2^2+2\sigma_{12}^2$ and $\E(Y_1^2Y_2)=\E(Y_1Y_2^2)=0$.
\end{theorem}

An application of \citeauthor{Isserlis1918}'s Theorem \ref{theo.Isserlis} gives

\begin{corollary}
\label{cor.N2}
Let $(X_1,X_2)\sim N_2\left(\left({\mu_1\atop\mu_2}\right),\left({\sigma_1^2\atop\sigma_{12}}~{\sigma_{12}\atop\sigma_2^2}\right)\right)$. Then, the covariance of $X_1^2$ and $X_2^2$ is $\Cov(X_1^2,X_2^2)=2\sigma_{12}\left(\sigma_{12}+2\mu_1\mu_2\right)$.
\end{corollary}

\section{\ \ Proofs}
\label{append.pr}

\begin{pr}{Proof of Theorem \ref{theo.bias,var}}
(a) By definition of $\ejl$,
\[
\E(\hmuos)=\frac{1}{n}\sum_{\mathclap{j=1}}^{k}\frac{1}{n-n_j}\sum_{\mathclap{l\ne j}}n_l\sum_{\mathclap{i\in S_j}}\E(\LL_{i,l})
          =\frac{1}{n}\sum_{\mathclap{j=1}}^{k}\frac{n_j}{n-n_j}\sum_{\mathclap{l\ne j}}n_l\ejl.
\]
The result arises by $n_j/n=p_j$ for all $j=1,\ldots,k$.
\medskip

\noindent
(b) Write $\Sigma_j=\sum_{l\ne j}n_l\sum_{i\in S_j}\LL_{i,l}$ and $\Sigma_{j,l}=\sum_{i\in S_j}\LL_{i,l}$. Hence, $\hmuos=\frac{1}{n}\sum_{j=1}^{k}\frac{1}{n-n_j}\Sigma_j$ and $\Sigma_j=\sum_{l\ne j}n_l\Sigma_{j,l}$. By straightforward calculations
\begin{equation}
\label{eq.var(hmuos)}
\begin{split}
\Var(\hmuos)=\frac{1}{n^2}&\left\{\sum_{\mathclap{j=1}}^k \frac{\Var(\Sigma_j)}{(n-n_j)^2}+\mathop{\sum\sum}\limits_{\mathclap{j\ne j'}}\,\,\frac{\Cov(\Sigma_j,\Sigma_{j'})}{(n-n_j)(n-n_{j'})}\right\}.
\end{split}
\end{equation}
The variance of $\Sigma_j$ is $\Var(\Sigma_j)=\sum_{l\ne j} n_l^2\Var(\Sigma_{j,l})+\sum\sum_{l\ne l'\colon l,l'\ne j}n_l n_{l'}\Cov(\Sigma_{j,l},\Sigma_{j,l'})$. We compute $\Var(\Sigma_{j,l})=\sum_{i\in S_j}\Var(\LL_{i,l})+\sum\sum_{i,i'\in S_j\colon i\ne i'}\Cov(\LL_{i,l},\LL_{i',l})=n_j\sfv_{j;l}+n_j(n_j-1)\sfc_{j;l}$, see \ref{eq.exch1} and \ref{eq.exch2}. Also, we compute the covariance of $\Sigma_{j,l}$ and $\Sigma_{j,l'}$, $\Cov(\Sigma_{j,l},\Sigma_{j,l'})=\sum_{i\in S_j}\Cov(\LL_{i,l},\LL_{i,l'})+\sum\sum_{i,i'\in S_j\colon i\ne i'}\Cov(\LL_{i,l},\LL_{i',l'})=n_j\sfc_{j;l,l'}$, see \ref{eq.exch3} and \ref{eq.exch5}. Thus,
\begin{equation}
\label{eq.var(Sj)}
\Var(\Sigma_j)=\sum_{\mathclap{l\ne j}} n_l^2\big(n_j\sfv_{j;l}+n_j(n_j-1)\sfc_{j;l}\big)+\mathop{\sum\sum}\limits_{\mathclap{l\ne l'\colon l,l'\ne j}}n_jn_l n_{l'}\sfc_{j;l,l'}.
\end{equation}
The covariance of $\Sigma_j$ and $\Sigma_{j'}$ is $\Cov(\Sigma_j,\Sigma_{j'})=\sum_{l\ne j}\sum_{l'\ne j'}n_l n_{l'}\Cov(\Sigma_{j,l},\Sigma_{j',l'})=\sum_{l\ne j,j'}n_l^2\Cov(\Sigma_{j,l},\Sigma_{j',l})+\sum\sum_{l\ne j,l'\ne j'\colon l\ne l'}n_l n_{l'}\Cov(\Sigma_{j,l},\Sigma_{j',l'})$. Now we compute $\Cov(\Sigma_{j,l},\Sigma_{j',l})=\sum_{i\in S_j}\sum_{i'\in S_{j'}}\Cov(\LL_{i,l},\LL_{i',l})=n_jn_{j'}\sfc_{j.j';l}$. For $\Cov(\Sigma_{j,l},\Sigma_{j',l'})$ when $l\ne l'$ we distinguish the following cases: If $l\ne j'$ and $l'\ne j$, $\Cov(\Sigma_{j,l},\Sigma_{j',l'})=0$, see \ref{eq.exch5}; if $l\ne j'$ and $l'=j$, $\Cov(\Sigma_{j,l},\Sigma_{j',j})=\sum_{i\in S_j}\sum_{i'\in S_{j'}}\Cov(\LL_{i,l},\LL_{i',j})=n_jn_{j'}\sfc_{j,j';l,j}$, see \ref{eq.exch4}; and if $l=j'$ and $l'\ne j$, similarly, $\Cov(\Sigma_{j,j'},\Sigma_{j',l'})=n_jn_{j'}\sfc_{j,j';j',l'}$. Thus, for each $j\ne j'$
\begin{equation}
\label{eq.cov(Sj,Sj')}
\Cov(\Sigma_j,\Sigma_{j'})
=\sum_{\mathclap{l\ne j,j'}}n_ln_jn_{j'}(n_l\sfc_{j,j';l}+2n_j\sfc_{j,j';l,j})
\end{equation}

Combining \eqref{eq.var(hmuos)}--\eqref{eq.cov(Sj,Sj')},
\[
\begin{split}
\Var(\hmuos)=\frac{1}{n^2}&\left\{\sum_{\mathclap{j=1}}^k \frac{\sum\limits_{\mathclap{l\ne j}} n_l^2\big(n_j\sfv_{j;l}+n_j(n_j-1)\sfc_{j;l}\big)+\mathop{\sum\sum}\limits_{\mathclap{l\ne l'\colon l,l'\ne j}}n_jn_l n_{l'}\sfc_{j;l,l'}}{(n-n_j)^2}\right.\\
                          &\left.~~+\mathop{\sum\sum}\limits_{\mathclap{j \ne j'}}\frac{\,\,\,\,\,\sum\limits_{\mathclap{l\ne j,j'}}n_ln_jn_{j'}(n_l\sfc_{j,j';l}+2n_j\sfc_{j,j';l,j})}{(n-n_j)(n-n_{j'})}\right\}.
\end{split}
\]
Using $n_j=np_j$ and $od_j=p_j/(1-p_j)$ for all $j$, after some algebra, the proof is completed.
\end{pr}

\begin{pr}{Proof of Lemma \ref{lem.moments}}
Since the joint distribution function of $X_j$s is unknown, the only information that we have is with respect to the moments of $X_j$s. Therefore, the only forms of estimators that we know to have expected values equal to linear combinations of $\mu^2$, $\sigma^2$ and $\sfc$ are the linear combination of $X_j^2$s and $X_jX_{j'}$s. Consider $\delta=\sum_{j=1}^n a_jX_j^2 + \sum\sum_{1\le j<j'\le n} b_{j,j'}X_jX_{j'}$. Observe that $\E(X_j^2)=\sigma^2+\mu^2$ and $\E(X_jX_{j'})=\sfc+\mu^2$. Thus, setting $a=\sum_{j=1}^n a_j$ and $b=\sum\sum_{1\le j<j'\le n} b_{j,j'}$, we get
\[
\E(\delta)=a(\sigma^2+\mu^2)+b(\sfc+\mu^2),
\]
completing the proof.
\end{pr}

\begin{pr}{Proof of Corollary \ref{cor.moments}}
(a)
We compute $\Var(\overline{X})=\frac{1}{n}\sigma^2+\frac{n-1}{n}\sfc$. Let now $\delta$ be an unbiased estimator of the variance of $\overline{X}$, $\delta=\sum_{j=1}^n a_jX_j^2 + \sum\sum_{1\le j<j'\le n} b_{j,j'}X_jX_{j'}$. Setting $a=\sum_{j=1}^n a_j$ and $b=\sum\sum_{1\le j<j'\le n} b_{j,j'}$, we get
\[
a(\sigma^2+\mu^2)+b(\sfc+\mu^2)=\frac{1}{n}\sigma^2+\frac{n-1}{n}\sfc
\Rightarrow
\left\{a=\frac{1}{n},\ b=1-\frac{1}{n},\ a+b=0\right\},
\]
a contradiction.
\medskip

\noindent
(b) Using the same arguments as in (a), the proof is completed.
\end{pr}

\begin{pr}{Proof of Theorem \ref{theo.consistency}}
Let $Z\sim F_j$ be independent of the elements of the $l$th source. Using H\"{o}lder inequality and \ref{cond2}, $\sfv_{j,l}\le\E[L(Z,\hd_{l,n})]^2\le\E\{[L(Z,\hd_{l,n})]^{2+\theta}\}^{2/(2+\theta)}\le M^{2/(2+\theta)}=M^*$ for all $j\ne l$. Thus, $\sum_{j=1}^k \frac{od_j^2}{np_j}(\sum_{l\ne j} p_l^2\sfv_{j;l})\le\frac{M^*}{n}\sum_{j=1}^k \frac{od_j^2}{p_j}(\sum_{l\ne j} p_l^2)=O(1/n)$. An application of Cauchy--Schwarz inequality gives $|\sfc_{j,l,l'}|\le\sfv_{j,l}\le M^*$. So, $\sum_{{j=1}}^k \frac{od_j^2}{np_j}({\sum\sum}_{{l\ne l'\colon l,l'\ne j}}p_lp_{l'}\sfc_{j;l,l'})\le\frac{M^*}{n}\sum_{{j=1}}^k \frac{od_j^2}{p_j}({\sum\sum}_{l\ne l'\colon l,l'\ne j}p_lp_{l'})=O(1/n)$.

It is remains to prove that $\sfc_{j;l}$, $\sfc_{j,j',l}$, $\sfc_{j,j';,l,j}\to0$ as $n\to\infty$. Let $Z,Z'$ are iid from $F_j$ and are independent of the elements of the $l$th source. Consider the sequence of random vectors $(Z,Z',\hd_{l,n})$ with respect to $n$. Then, \ref{cond1} gives $(Z,Z',\hd_{l,n})\rightsquigarrow(Z,Z',d_l)$ as $n\to\infty$. Since $L$ is continuous, the maps $(Z,Z',\hd_{l,n})\mapsto L(Z,\hd_{l,n})$ and $(Z,Z',\hd_{l,n})\mapsto L(Z,\hd_{l,n})L(Z',\hd_{l,n})$ are continuous. Using the Continuous Mapping Theorem, $L(Z,\hd_{l,n})\rightsquigarrow L(Z,d_l)$ and $L(Z,\hd_{l,n})L(Z',\hd_{l,n})\rightsquigarrow L(Z,d_l)L(Z',d_l)$ as $n\to\infty$. Observe that $\E|L(Z,\hd_{l,n})|^{1+(1+\theta)}\le M$, see \ref{cond2}, so the sequence $L(Z,\hd_{l,n})$ is uniformly integrable. Hence, $\E|L(Z,d_l)|<\infty$ and $\E[L(Z,\hd_{l,n})]\to\E[L(Z,d_l)]$ as $n\to\infty$ \citep[see, e.g.,][p.~338]{Billingsley1995}. Similarly, $\E|L(Z',d_l)|<\infty$ and $\E[L(Z',\hd_{l,n})]\to\E[L(Z',d_l)]$ as $n\to\infty$. Using Cauchy--Schwarz inequality we obtain $\E|L(Z,\hd_{l,n})L(Z',\hd_{l,n})|^{1+\theta/2}\le(\E|L(Z,\hd_{l,n})|^{2+\theta})^{1/2}(\E|L(Z',\hd_{l,n})|^{2+\theta})^{1/2}\le M$. So, $L(Z,\hd_{l,n})L(Z',\hd_{l,n})$ is uniformly integrable. Therefore, $\E|L(Z,d_l)L(Z',d_l)|<\infty$ and $\E[L(Z,\hd_{l,n})L(Z',\hd_{l,n})]\to\E[L(Z,d_l)L(Z',d_l)]$ as $n\to\infty$. Moreover, $L(Z,d_l)$ and $L(Z',d_l)$ are independent. So, $\E[L(Z,\hd_{l,n})L(Z',\hd_{l,n})]\to\E[L(Z,d_l)]\E[L(Z',d_l)]$ as $n\to\infty$. From the preceding analysis we have that $\sfc_{j;l}\to0$ as $n\to\infty$. Using the same arguments as above it follows that $\sfc_{j,j',l}$, $\sfc_{j,j';,l,j}\to0$ as $n\to\infty$, and the proof is completed.
\end{pr}

\begin{pr}{Proof of Equations \eqref{eq.normal,sq-abs1}--\eqref{eq.normal,sq-abs3}}
Let $j,j',l,l'$ are four distinct indices. Assume that $Z_{(j)},Z'_{(j)}\in \CZ_{S_j}$ (with $Z_{(j)}\ne Z'_{(j)}$) and $Z_{(j')}\in \CZ_{S_{j'}}$. Consider the following random vectors $(X_1,X_2)=(Z_{(j)}-\overline{Z}_l,Z'_{(j)}-\overline{Z}_l)\sim N_2\left(\left({\mu_j-\mu_l\atop\mu_j-\mu_l}\right), \left({\sigma_j^2+\sigma_l^2/n_l\atop\sigma_l^2/n_l}~{\sigma_l^2/n_l\atop\sigma_j^2+\sigma_l^2/n_l}\right)\right)$, $(X_3,X_4)=(Z_{(j)}-\overline{Z}_l,Z_{(j)}-\overline{Z}_{l'})\sim N_2\left(\left({\mu_j-\mu_l\atop\mu_j-\mu_{l'}}\right), \left({\sigma_j^2+\sigma_l^2/n_l\atop\sigma_j^2}~{\sigma_j^2\atop\sigma_j^2+\sigma_{l'}^2/n_{l'}}\right)\right)$, $(X_5,X_6)=(Z_{(j)}-\overline{Z}_l,Z_{(j')}-\overline{Z}_l)\sim N_2\left(\left({\mu_j-\mu_l\atop\mu_{j'}-\mu_l}\right), \left({\sigma_j^2+\sigma_l^2/n_l\atop\sigma_l^2/n_l}~{\sigma_l^2/n_l\atop\sigma_{j'}^2+\sigma_l^2/n_l}\right)\right)$, $(X_7,X_8)=(Z_{(j)}-\overline{Z}_l,Z_{(j')}-\overline{Z}_j)\sim N_2\left(\left({\mu_j-\mu_l\atop\mu_{j'}-\mu_j}\right), \left({\sigma_j^2+\sigma_l^2/n_l\atop-\sigma_j^2/n_j}~{-\sigma_j^2/n_j\atop\sigma_{j'}^2+\sigma_j^2/n_j}\right)\right)$. If $X\sim N(\mu,\sigma^2)$, then $\Var(X^2)=2\sigma^2(\sigma^2+2\mu^2)$. Since $\sfv_{j,l}=\Var(X_1^2)$, \eqref{eq.normal,sq-abs2(a)} follows.  Observe that $\sfc_{j;l}=\Cov(X_1^2,X_2^2)$, $\sfc_{j;l,l'}=\Cov(X_3^2,X_4^2)$, $\sfc_{j,j';l}=\Cov(X_5^2,X_6^2)$ and $\sfc_{j,j';l,j}=\Cov(X_7^2,X_8^2)$. Hence, an application of Corollary \ref{cor.N2} proves \eqref{eq.normal,sq-abs2(b)}. Finally, If $X\sim N(\mu,\sigma^2)$, then
$\E|X|=\mu[1-2\Phi(-\mu/\sigma)]+\sigma(2/\pi)^{1/2}\exp\{-\mu^2/(2\sigma^2)\}$. Because $\ejl=\E|X_1|=\mu_{j;l}\left[1-2\Phi\left(-{\mu_{j;l}}/{\sigma_{j;l}}\right)\right]+\sigma_{j;l}\sqrt{2/\pi}\exp\left(-{\mu_{j;l}^2}/{2\sigma_{j;l}^2}\right)$, using \eqref{eq.OOSer}, \eqref{eq.normal,sq-abs3} follows, completing the proof.
\end{pr}

\begin{pr}{Proof of Corollary \ref{cor.N2}}
Consider $Y_1=X_1-\mu_1$ and $Y_2=X_2-\mu_2$. Then, $(Y_1,Y_2)\sim N_2\left(\bm{0},\left({\sigma_1^2\atop\sigma_{12}}~{\sigma_{12}\atop\sigma_2^2}\right)\right)$. Thus, $\Cov(X_1^2,X_2^2)=\Cov(Y_1^2+2\mu_1Y_1+\mu_1^2,Y_2^2+2\mu_2Y_2+\mu_2^2)=\Cov(Y_1^2,Y_2^2)+2\mu_2\Cov(Y_1^2,Y_2)+2\mu_1\Cov(Y_1,Y_2^2)+4\mu_1\mu_2\Cov(Y_1,Y_2)$. A simple application of \citeauthor{Isserlis1918}'s Theorem completes the proof.
\end{pr}

\end{appendices}


\begin{thebibliography}{}
\bibitem[Afendras and Markatou, 2016]{AM2016}
    Afendras, G.\ and  Markatou, M. (2016).
    Optimality of Training/Test Size and Resampling Effectiveness of Cross-Validation Estimators of the Generalization Error.
    \textit{arXiv}:1511.02980v1 [math.ST]

\bibitem[Arlot and Celisse, 2010]{AC2010}
    Arlot, S.\ and Celisse, A. (2010).
    A survey of cross-validation procedures for model selection.
    \textit{Statistics Surveys},
    \textbf{4}, 40--79.

\bibitem[Ben-David et al., 2010]{B-DBCKPV2010}
    Ben-David, S., Blitzer, J, Crammer, K., Kulesza, A., Pereira, F.\ and Vaughan, J.W. (2010).
    A theory of learning from different domains.
    \textit{Mach.\ Learn.},
    \textbf{79}, 151--175.

\bibitem[Bengio and Grandvalet, 2004]{BG2004}
    Bengio, Y.\ and Grandvalet, Y. (2004).
    No unbiased estimator of the variance of $k$-fold cross validation.
    \textit{Journal of Machine Learning Research},
    \textbf{5}, 1089--1105.

\bibitem[Breiman et al., 1984]{BFOS1984}
    Breiman, L., Friedman, J.H., Olshen, R.A.\ and Stone, C.J. (1984).
    \textit{Classification and Regression Trees}.
     Wadsworth, Belmont, California.

\bibitem[Billingsley, 1995]{Billingsley1995}
    Billingsley, P. (1995).
    {Probability and Measure.}
    Third Edition.
    Wiley Series in Probability and Mathematical Statistics.

\bibitem[Geisser, 1975]{Geisser1975}
    Geisser, S. (1975).
    The predictive sample reuse method with applications.
    \textit{Journal of the American Statistical Association},
    \textbf{70}(350), 320--328.

\bibitem[Geras and Sutton, 2013]{GS2013}
    Geras, K.\ and Sutton, C. (2013).
    Multiple-source cross-validation.
    \textit{Proceedings of the 30 th International Conference on Machine Learning},
    Atlanta, Georgia, USA, 2013.
    JMLR: W\&CP, \textbf{28}(3), 1292--1300.

\bibitem[Isserlis, 1918]{Isserlis1918}
    Isserlis, L. (1918).
    On a formula for the product-moment coefficient of any order of a normal frequency distribution in any number of variables.
    \textit{Biometrika},
    \textbf{12}, 134--139.

\bibitem[Markatou et al., 2005]{MTBH2005}
    Markatou, M., Tian, H, Biswas, S.\ and Hripcsak, G. (2005).
    Analysis of Variance of Cross-Validation Estimators of the Generalization Error.
    \textit{J.\ Mach.\ Learn.\ Res.},
    \textbf{6}, 1127--1168.

\bibitem[Nadeau and Bengio, 2003]{NB2003}
    Nadeau, C.\ and Bengio, Y. (2003).
    Inference for the generalization error.
    \textit{Mach.\ Learn.},
    \textbf{52}, 239--281.

\bibitem[Stone, 1974]{Stone1974}
    Stone, M. (1974).
    Cross-validatory choice and assessment of statistical predictions.
    \textit{Journal of the Royal Statistical Society. Series B},
    \textbf{36}(2), 111--147.

\bibitem[Stone, 1977]{Stone1977}
    Stone, M. (1977).
    Asymptotics for and against cross-validation.
    \textit{Biometrika},
    \textbf{64}(1), 29--35.


\end{thebibliography}
\end{document}